\documentclass[10pt]{article}
\usepackage{amsmath,amssymb}

\usepackage{graphicx}

\usepackage{amsmath}
\usepackage{tabularx} 
\usepackage{setspace}  
 
\newcommand{\ba}{\begin{align}} 
\newcommand{\ea}{\end{align}}

\usepackage{amsfonts} 
\usepackage{amssymb}

\def\be{\begin{eqnarray}}
\def\ee{\end{eqnarray}}
\def\bee{\begin{eqnarray*}}
\def\eee{\end{eqnarray*}}

\newtheorem{thm}{Theorem}

\newtheorem{lemma}[thm]{Lemma}

\newtheorem{defn}[thm]{Definition}

\def\la{\langle}
\def\ra{\rangle}
\def\E{\mathbb{E}}
\def\P{\mathbb{P}}

\def\bmx{\begin{pmatrix}}
\def\emx{\end{pmatrix}}

\newcounter{exercise}
\title{The fluid limit of a random graph model for a shared ledger}

    \author{Christopher King \\
    \\
  {\small      Department of Mathematics} \\
  {\small      Northeastern University} \\
  {\small      Boston MA 02115}
}

\begin{document}

     \maketitle

     \begin{abstract}
A shared ledger is a record of transactions that can be updated by any member of a
group of users. The notion of independent and consistent record-keeping in a  shared ledger  is important for blockchain and
more generally for distributed ledger technologies. In this paper we analyze the growth of a model for the  tangle,
which is the shared ledger protocol used as the basis for the IOTA cryptocurrency.
The model is a random directed acyclic graph, and its growth is described by a non-Markovian stochastic process.
We derive a delay differential equation for the fluid model which describes the tangle at high arrival rate.
We prove convergence in probability of the tangle process to the fluid model,
and also prove global stability of the fluid model. The convergence proof relies on martingale techniques.
     \end{abstract}

\section{Introduction}
In this paper we analyze a stochastic growth process for a family of directed acyclic graphs, and show that the fluid limit of this
process is described by a delay differential equation. This stochastic process describes a type of shared ledger which
was introduced as the foundation of the cryptocurrency IOTA \cite{Popov}, and the result about the fluid limit was
used previously to analyze the persistence of competing transaction records in this ledger  \cite{FKS1}.
The main contribution of this paper is to provide a precise formulation of the results about the fluid limit.
We use martingale techniques to establish that the process converges weakly  in the limit where the arrival rate goes to infinity,
and show that the fluid limit is given as the solution of a suitable delay differential equation.
We also prove a convergence result for the  solutions of the delay differential equation.

\medskip
The term `shared ledger' refers to a record of transactions which may be amended independently by any member of a
group of users.
The goal of designing a shared ledger is to allow users to add transactions to the record without 
centralized control, while at the same time protecting the record against tampering by malicious agents. 
As a background to this topic 
we will review below the well-known blockchain protocol \cite{Nakamoto} which
involves linking blocks (collections of transactions) by complicated hash function computations.
If blocks are represented by vertices on a graph and
the hashing link between two blocks is represented by a directed edge between those vertices,
then the whole blockchain ledger can be viewed as a directed graph. This point of view leads to our random graph analysis,
and will form the basis for the stochastic process that will be analyzed in this paper.

\subsection{The blockchain protocol}
The blockchain technology underlying Bitcoin is a well-known implementation of a shared ledger which provides
security against malicious users \cite{Nakamoto}, \cite{Zheng}. Recall that the blockchain is
an ordered string of blocks, each containing several hundred transaction records; each block has a unique numerical ID (256 bits for Bitcoin)
that satisfies a challenging constraint.
The ID of a block is computed using a complicated hash function, and the input for the hash function
involves the block's own data, the ID of the previous block, and some extra bits which are chosen so that the output satisfies the constraint. 
Thus every block's ID depends on the data
of the previous block, and hence also on the data in all previous blocks. Therefore any change in the data of a 
block would change the IDs of all subsequent blocks, and the altered IDs would almost certainly not satisfy the 
tight constraint mentioned above. This failure would be a signal to all observers that the ledger had been altered,
and so the existence of a ledger with valid IDs for all blocks is its own guarantee of security.
The key mechanism for security is the difficulty of computing a valid ID for a block.
This task is called the proof of work, and requires finding an input to a complicated hash function which will produce
an output of the specified form.
The blockchain miners compete to find this inverse, and the first successful one adds the new block to the chain.

\medskip
One essential constraint in Bitcoin is that a new block can only be linked to the most recent block in the chain. This constraint
ensures that the blockchain is a linear graph. It also ensures that every transaction record in the chain is linked to all subsequent records, and indeed
the security of a transaction increases as later blocks are added (a typical rule of thumb is that a transaction record in blockchain
is `safe' after at least six subsequent blocks have been added). However this constraint leads to a `winner takes all' rule for the miners,
who must compete to be first to add a new block. Consequently mining has become a dedicated enterprise requiring
specialized technology, and there is much effort wasted (and energy expended) by the miners.

\subsection{Modifying the blockchain protocol}
There have been many proposed modifications of the blockchain protocol. In this paper we consider one such proposal \cite{Popov} which involves removing the
constraint that a new block can only be linked to the most recent block in the chain. Removing this constraint has several immediate consequences.
First,
there is no competition between miners, hence each user can perform their own proof of work (which is much easier than in Bitcoin) and there are no
rewards for adding a new block. Second, since
a new block can link to any previous block, the graph of links for the ledger is no longer linear,  and 
can be much more complicated than in Bitcoin. Furthermore,
since there are many possible ways to link a new block to the ledger, it is reasonable to view the ledger as a randomly growing
graph and to investigate its typical properties. We will pursue this point of view for the modification known as
the {\em tangle protocol} which was introduced in  \cite{Popov}.

\subsection{The tangle protocol}
In the tangle protocol \cite{Popov} each new block contains just one transaction. A new transaction
links to {\em two} existing transactions in the ledger (this is another change from the blockchain protocol), 
and the proof of work
uses the ID's of these two transactions as part of its input. Thus the ledger grows by the addition of transactions each with two directed edges
which link to existing transactions in the ledger. The resulting graph of links is a connected directed acyclic  graph (DAG).
The proof of work lasts for some amount of time $h$, so there is a delay between the time when a new transaction starts its proof of
work and the time when it is added to the DAG as a new vertex. This time $h$ is much shorter than 
in Bitcoin (where the proof of work on average lasts for ten minutes), however it plays a crucial role in the growth of the ledger.

\medskip
As in the blockchain, the security of a transaction increases as later transactions are added which are linked either directly or indirectly
to it. Although in principle a new transaction may choose to link to any existing transaction in the ledger,
it is advantageous to select two recently arrived transactions for linking. 
Transactions which have not yet been linked by subsequent transactions are called tips; 
in the tangle protocol all users select tips for linking.
Note that for blockchain the security relies on the community solving
one exceedingly difficult hash function inversion for each new block, 
whereas for the tangle the security relies on a large community of users each
performing relatively simple computations in parallel.

\subsection{Summary of results}
There are two time scales in this model, namely $\lambda^{-1}$, the average time between arrivals of new transactions,
and $h$, the duration of the proof of work. We are interested in how the average number of tips in the 
tangle depends on these quantities as the arrival rate $\lambda$ goes to infinity. Let $L(t)$ be the number of tips at time $t$. 
Several approaches to this question \cite{Popov}, \cite{FKS1} have shown that $L(t)$ is roughly
proportional to the product $\lambda h$, at least for large $\lambda$. We investigate in this paper the limit where $\lambda$ approaches
infinity, so we define the rescaled variable $B^{(\lambda)}(t) = \lambda^{-1} L(t)$. The first
result Theorem \ref{thm1} shows that $B^{(\lambda)}(t)$ converges in probability to a deterministic function $b(t)$
as $\lambda \rightarrow \infty$. The bound in Theorem \ref{thm1} also shows that the fluctuations in $|B^{(\lambda)}(t) - b(t)|$ are 
no larger than
$O(\lambda^{-1/2})$. The second result Theorem \ref{thm2} shows that the function $b(t)$ converges exponentially to
$2 h$ as $t \rightarrow \infty$. Putting these results together shows that for large $\lambda$ and large $t$,
$L(t)$ can be written approximately as $2 \lambda h + O(\lambda^{1/2})$.

\subsection{Outline of the paper}
In Section \ref{sec:DAG} we formulate a stochastic process for the number of tips on the DAG
which represents the tangle. In Section \ref{sec:fluid} we describe the fluid limit of the rescaled process
(the fluid limit refers to the limit where the arrival rate of new transactions goes to infinity)
and also describe how initial conditions can be consistently formulated for the process and the
fluid limit.
The main results of the paper Theorems \ref{thm1} and \ref{thm2} are stated in Section \ref{sec:results},
and are proved in Sections \ref{sec:pf} and \ref{sec:pf2}. Future directions of research
on this topic are discussed in Section \ref{sec:discuss}.

\section{Definition of the DAG model}\label{sec:DAG}
Let $G=(V,E)$ be a finite connected acyclic directed graph (DAG) where $V$ is the vertex set and $E$ is the edge set.
If an edge $e \in E$ is directed from vertex $x$ to vertex $y$ we will write $e = \la x, y \ra$ and say that $y$ is the head 
and $x$ is the tail. A {\em tangle} is a DAG with two additional properties: first,
there is a unique vertex which is not the tail in any edge -- this is called the genesis vertex.
Second, every vertex is the tail in at most two edges. 
The subset of vertices which are not the heads of any edges will be called the {\em tips} of the tangle.

\medskip
We will define a stochastic growth model for the tangle. The arrival rate of new transactions is denoted by $\lambda$, and
for simplicity we assume that new transactions are created
at the deterministic sequence of times $\{t_n = \lambda^{-1} n \,: \, n=1,2,\dots\}$.
At time $t_n$, two tips $x_1(n)$ and $x_2(n)$
on the tangle are selected for the proof of work by the new transaction (it is possible that $x_1(n) = x_2(n)$).
The proof of work 
lasts for a fixed length of time $h$. For simplicity we will 
assume that $\lambda$ is always chosen so that $\lambda h$ is an integer:
\be\label{def:h}
m = \lambda h
\ee
At time $t_n + h = t_{n+m}$ the new transaction is added to the tangle as a tip $y_n$,
and the two directed edges $\la y_n, x_1(n) \ra$ and $\la y_n, x_2(n) \ra$ are also added to the graph.
This is the only mechanism by which the tangle grows.

\medskip
Obviously the vertices $x_1(n)$ and $x_2(n)$ are no longer tips after time $t_n+h$,
however it may happen that these vertices had already ceased  to be tips at an earlier time, due to their being linked to some
other new transaction. We say that a tip is {\em pending}  if it has been selected for proof of work by
a transaction but has not yet been linked. We say that a tip is {\em free} if it is not pending.

\begin{defn}
\be
W_n &=& \mbox{number of pending tips at time $t_n$} \\ 
X_n &=& \mbox{number of free tips at time $t_n$} \\
L_n &=& W_n + X_n = \mbox{number of tips at time $t_n$} \\
U_n &=& \mbox{number of free tips selected for proof of work at time $t_n$}
\ee
\end{defn}

We have defined $U_n$ to be the number of the vertices $\{x_1(n), x_2(n)\}$
which are free at time $t_n$, so $U_n \in \{0,1,2\}$.  After selection these free vertices immediately become pending vertices, 
hence they will never contribute to any of the subsequent values $U_{n+1}, U_{n+2}, \dots$.
Furthermore at any time $n \ge m$ there are exactly $m$ new transactions
which are each in the process of carrying out their proof of work on two vertices on the graph
(this holds because $m = \lambda h$ and we assume that the value of $h$ is fixed and identical for all
users). Therefore the total number of pending vertices at any time $t_n$ (with $n \ge m$) 
is the sum of $\{U_n, U_{n-1}, \dots, U_{n-m+1}\}$, that is
\be\label{eqn-W}
W_n = \sum_{j=n-m+1}^n U_j \quad \mbox{for all $n \ge m$}
\ee
We also have the following evolution relations:
\be\label{update1}
X_{n+1} &=& X_n + 1 - U_{n+1} \nonumber \\
L_{n+1} &=& L_n + 1  - U_{n - m + 1}
\ee
and the relation
\be\label{L=WX}
L_n = W_n + X_n
\ee
We will discuss shortly how the processes $X_n$ and $L_n$ can be defined using appropriate initial conditions.
For this purpose it will be convenient to use the evolution relations (\ref{update1}) starting at $n=m$ and $n=2m$ respectively.
These relations imply that
\be\label{update1a}
X_{n+1} &=& X_m + \sum_{i=m}^n (1 - U_{i+1}) \quad \mbox{ for all $n \ge m$} \nonumber \\
L_{n+m+1} &=& L_{2 m} + \sum_{i=m}^n (1 - U_{i+1}) \quad \mbox{ for all $n \ge 2 m$}
\ee
Therefore 
\be\label{update1b}
L_{n+m+1} = X_{n+1} + L_{2 m} - X_m
\ee
and we have from (\ref{L=WX}), (\ref{eqn-W}) and (\ref{update1a}) that
\be
L_{2 m} = W_{2 m} + X_{2 m} = \sum_{i=m+1}^{2 m} U_i + X_m + \sum_{i=m}^{2 m-1} (1 - U_{i+1})
= X_m + m
\ee
So we conclude form (\ref{update1b}) that for all $n \ge 2 m$,
\be\label{LXm}
L_{n} = X_{n-m} + m
\ee
We will refer to $(W_n,X_n,L_n)$ as the tangle process. 
It follows from (\ref{LXm}) that the process $\{L_n \}$ (for $n \ge 2 m$) is in fact fully determined by $\{X_n : n \ge m\}$.

\subsection{The random tip growth model}
The remaining ingredient in the definition of the DAG model is the method of choosing vertices $x_1(n)$ and $x_2(n)$.
We will assume in this paper that the tips $x_1(n)$ and $x_2(n)$ are chosen independently and uniformly from the set of tips,
and we call this the {\em random tip growth} (RTG) model for the tangle process. 
Thus the numbers $\{U_n\}$ are random variables whose distributions depend on the number of tips at time $t_n$.
The RTG model is one of the tip selection algorithms discussed in \cite{Popov}, \cite{FKS2}, and it is expected that the 
fluid limit results presented in this paper can be extended to those other tip selection algorithms.

\medskip
We denote by ${\cal F}(\lambda,n)$ the $\sigma$-algebra generated by $\{U_1,\dots,U_n\}$:
\be
{\cal F}(\lambda,n) = \sigma \left(U_1,\dots,U_n \right)
\ee
It follows from (\ref{update1a})  that $X_n - X_0$ and $L_{n+m} - L_m$ are measurable with respect to ${\cal F}(\lambda,n)$. 
We also have the filtration relation
\be
{\cal F}(\lambda,n_1) \subset {\cal F}(\lambda,n_2) \quad \mbox{for all $n_1 < n_2$}
\ee
The conditional distribution of the random variable $U_{n+1}$ for the RTG model is
\be\label{dist-U}
\P(U_{n+1}=2\,|\, {\cal F}(\lambda,n)) &=& \frac{X_n (X_n - 1)}{{L_n}^2} \nonumber \\
\P(U_{n+1}=0\,|\, {\cal F}(\lambda,n)) &=& \frac{{W_n}^2}{{L_n}^2} \nonumber \\
\P(U_{n+1}=1\,|\, {\cal F}(\lambda,n)) &=& \frac{2 W_n X_n + X_n}{{L_n}^2}
\ee
Note that
\be\label{exp-U}
\E[U_{n+1} \,|\, {\cal F}(\lambda,n)] = 2 \, \frac{X_n}{L_n} - \frac{X_n}{L_n^2}
\ee
It is clear that $(X_n,L_n)$ is not a Markov process, as the distribution of
$L_{n+1}$ depends on $U_{n-m+1}$, which in turn depends on $(X_{n-m}, L_{n-m})$ through (\ref{dist-U}).

\subsection{Generating the process from initial conditions}\label{sec:gen-init-stoch}
The stochastic process $(X_n,L_n)$ defined by (\ref{update1}) and (\ref{dist-U}) must be supplemented with initial conditions 
in order to be well-defined. This is done most easily by assigning values to the variables 
$(U_1,\dots,U_m)$ and $X_m$. Once these assignments have been made, the distribution of the process $(X_n, L_n)$
is determined for all $n \ge m$, as will be explained below. In particular the
variables $X_0,\dots,X_{m-1}$ and $L_0,\dots,L_{m-1}$
do not play any role, and we will ignore their values.

\medskip
Let $(u_1,\dots,u_m)$ be a sequence with $u_i \in \{0,1,2\}$ for all $i=1,\dots,m$, and let
$\xi_m \ge 1$ be an integer. Then we assign as initial conditions
\be\label{init-proc1}
U_i &=& u_i, \quad i=1,\dots,m \nonumber \\
W_m &=& \sum_{i=1}^m u_i \nonumber \\
X_m &=& \xi_m \nonumber \\
L_m &=& X_m + W_m = \xi_m + \sum_{i=1}^m u_i
\ee
To see how this generates the process for $n > m$, note from (\ref{dist-U}) that the distribution of $U_{m+1}$ is determined
by (\ref{init-proc1}), and is well-defined since $L_m \ge \xi_m \ge 1$. Given $U_{m+1}$, we compute
\be\label{init-proc2}
W_{m+1} &=& \sum_{i=2}^{m+1} U_i\nonumber  \\
X_{m+1} &=& X_m + 1 - U_{m+1} \nonumber \\
L_{m+1} &=& W_{m+1} + X_{m+1}
\ee
We have
 $L_{m+1} \ge X_{m+1} \ge X_m -1$, and
also the formulas (\ref{dist-U}) imply that if $X_m =1$ then $U_{m+1} \in \{0,1\}$, and thus $L_{m+1} \ge X_{m+1} \ge X_m$.
Therefore $L_{m+1} \ge 1$, and so the distribution of $U_{m+2}$ is well-defined, and this random variable can be used to
compute $W_{m+2}, X_{m+2}, L_{m+2}$.
This bootstrap
continues and allows us to generate the whole process starting from the initial conditions (\ref{init-proc1}), 
and furthermore the resulting process satisfies $X_n \ge 1$ for all $n \ge m$.
Hence from
(\ref{LXm}) we also deduce 
that $L_n \ge m + 1$ for all $n \ge 2 m$.

\medskip
From (\ref{init-proc2}) and (\ref{init-proc1}) we also deduce that
\be
L_{m+1} = \xi_m + 1 + \sum_{i=2}^m u_i
\ee
and so $L_{m+1}$ is also fixed by the initial conditions. The same is true for $L_2,\dots,L_m$, and we have the formula
\be\label{init-proc3}
L_{m+j} = \xi_m + j + \sum_{i=j+1}^m u_i \quad \mbox{for $j=0,\dots,m$}
\ee

\section{The fluid limit}\label{sec:fluid}
Given the process $\{X_n, L_n\}$ we rescale variables and define for all $t > 0$
\be\label{rescale1}
A^{(\lambda)}(t) = \lambda^{-1} \, X_{n(t)}, \qquad
B^{(\lambda)}(t) = \lambda^{-1} \, L_{n(t)}, \qquad 
\mbox{where $n(t) = \lfloor \lambda \, t \rfloor$}
\ee
The variables $(A^{(\lambda)}(t), B^{(\lambda)}(t))$ are piecewise constant in the intervals $[t_n,t_{n+1})$,
and change by at most $\pm \lambda^{-1}$ at each time $t_n$. Thus it is reasonable that in the limit
$\lambda \rightarrow \infty$ these variables will converge to continuous functions $a(t)$ and $b(t)$.
Furthermore after rescaling (\ref{update1}) the evolution equations become
\be\label{rescale2}
\frac{A^{(\lambda)}(t_{n+1}) - A^{(\lambda)}(t_n)}{t_{n+1} - t_n} = 1 - U_{n+1} \nonumber \\
\frac{B^{(\lambda)}(t_{n+1}) - B^{(\lambda)}(t_n)}{t_{n+1} - t_n} = 1 - U_{n -m+1}
\ee
The left sides of (\ref{rescale2}) are expected to converge to $a'(t)$ and $b'(t)$ as $\lambda \rightarrow \infty$,
so it is reasonable to expect that the fast variations on the right side will be averaged out in the limit,
leaving the expected values of the variables $U_{n+1}$ and $U_{n-m+1}$.
From (\ref{exp-U}) we have
\be\label{rescale3}
\E[U_{n+1} \,|\, {\cal F}(\lambda,n)] = 
2 \, \frac{A^{(\lambda)}(t_n)}{B^{(\lambda)}(t_n)} - \lambda^{-1} \, \frac{A^{(\lambda)}(t_n)}{B^{(\lambda)}(t_n)} 
 \simeq  2 \, \frac{A^{(\lambda)}(t_n)}{B^{(\lambda)}(t_n)} 
 \simeq  2 \, \frac{a(t)}{b(t)}
\ee
and similarly
\be\label{rescale4}
\E[U_{n-m+1} \,|\, {\cal F}(\lambda,n)] 
 \simeq  2 \, \frac{A^{(\lambda)}(t_n - h )}{B^{(\lambda)}(t_n - h)} 
 \simeq  2 \, \frac{a(t - h)}{b(t - h)}
\ee
Assuming that the right sides of (\ref{rescale2}) converge to these average values, we are led to the following 
pair of coupled delay differential equations for the fluid limit:
\be\label{DDE1}
\frac{d a}{d t} = 1 - 2\, \frac{a(t)}{b(t)}, \qquad \frac{d b}{d t} = 1 - 2\, \frac{a(t - h)}{b(t - h)}
\ee

\subsection{Delay differential equations}
The equations (\ref{DDE1}) must be supplemented with suitable initial conditions. We will say that the combination 
$\alpha = (a(h), \{u(t) \,:\, 0 \le t \le h\})$ 
is a {\em DDE initial condition} if $a(h) > 0$, $u(t)$ is integrable and 
\be\label{DDE-init1}
0 \le u(t) \le 2 \quad \mbox{for all $0 \le t \le h$}
\ee
These initial conditions can be used to define a solution of the fluid equations (\ref{DDE1}) for $t \ge 2h$,
in the same way as the initial conditions (\ref{init-proc1}) were used to construct the tangle process.
The idea is that the function $u(t)$ plays the same role as the initial sequence $\{u_i\}$ for the discrete process.
Thus we first define the initial value $b(h)$ as
\be\label{init-b0}
b(h) = a(h) + \int_{0}^h u(s) \, d s,
\ee
and we then define $b(t)$ for $h \le t \le 2h$ as the solution of the delay equation
\be
\frac{d b}{d t} = 1 - u(t-h)
\ee
This leads to the solution
\be\label{init-b1}
b(t) &=& b(h) + t - h  - \int_{0}^{t-h} u(s) \, d s \nonumber \\
&=& a(h) + t - h + \int_{t-h}^h u(s) \, d s \quad  \mbox{for $h \le t \le 2h$}
\ee
We then compute $a(t)$ for $t \in [h,2h]$ as the solution of the equation
\be
\frac{d a}{d t} = 1 - 2 \frac{a(t)}{b(t)}
\ee
which gives
\be\label{init-a1}
a(t) = P(h,t)^{-1} a(h) + P(h,t)^{-1} \, \int_{h}^t P(h,s) \, d s, \quad \mbox{for $t \in [h,2h]$}
\ee
where
\be\label{def:P}
P(x,y) = \exp \left(2 \int_{x}^y b(s)^{-1} \, d s\right)
\ee
Note that (\ref{init-b1}) implies $b(t) \ge a(h) > 0$ for all $h \le t \le 2h$, so (\ref{def:P}) is well-defined for $(x,y) = (h,t)$ with $t$ in this interval,
and (\ref{init-a1}) also implies that $a(t) > 0$ for all $h \le t \le 2h$.
The equation (\ref{init-b1}) also implies that $b(2h) = a(h) + h$.
Having obtained the functions $(a(t),b(t))$ in the interval $[h,2h]$, we then extend the solutions to the interval
$[2h,3h]$ by first defining
\be\label{init-b2}
b(t) = a(t-h) + h \quad \mbox{for all $2h \le t \le 3h$}
\ee
and then solving the differential equation for $a(t)$ to obtain
\be\label{init-a2}
a(t) = P(2h,t)^{-1} a(2h) + P(2h,t)^{-1} \, \int_{2h}^t P(2h,s) \, d s, \quad \mbox{for $2h \le t \le 3h$}
\ee
From (\ref{init-b2}) we have $b(t) \ge h$, and thus $P(2h,t)$ is well-defined for $2h \le t \le 3h$, and again
implies positivity of $a(t)$.  This construction can be continued in the same way for subsequent intervals $[3h,4h],\dots$,
and produces a solution of the equations (\ref{DDE1}) for all $t > 2h$.
We collect together our results about this solution in the following lemma.

\begin{lemma}\label{lem1}
Let $\alpha$ be a DDE initial condition. There are unique functions
$(a(t),b(t))$ defined for all $t > h$ which satisfy the equations (\ref{init-b1}) and (\ref{init-a1}) in the interval $[h,2h]$, 
and which satisfy the differential equations (\ref{DDE1}) for all $t > 2 h$. For $t \ge 2h$ the solutions
also satisfy the following conditions:
\be\label{lem1:eqn1}
1) && a(t) \ge 0, \\
2) && b(t) = h + a(t-h) \\
3) && b(t) \ge h \\
4) && b(t) - a(t) = \int_{t-h}^t 2 \frac{a(s)}{b(s)} \, d s \\
5) && 0 \le b(t) - a(t) \le 2h
\ee
\end{lemma}

\medskip
\noindent{\em Proof of Lemma \ref{lem1}}:
the formulas (\ref{init-b1}) and (\ref{init-a1}) show that $(a(t),b(t))$ is uniquely defined and differentiable in the interval $(h,2h)$,
and is continuous at $t=2h$. The iterative construction outlined above produces a unique differentiable solution in every interval
$(j h, (j+1) h)$ for $j=2,3,\dots$. The solution is clearly continuous at $t= j h$ for all $j \ge 2$.
It is also differentiable at $t= j h$ for all $j \ge 3$ because it satisfies the differential equations (\ref{DDE1}) in both intervals
$((j-1) h, jh)$ and $(j h, (j+1) h)$. Properties (1), (2), (3) follow by construction. To see that Property (4) holds, let $c(t) = b(t) - a(t)$ and
consider first the interval $[h,2h]$, where we have
\be
c'(t) = b'(t) - a'(t) = 2 \frac{a(t)}{b(t)} - u(t-h)
\ee
Therefore for some constant $K$ we have
\be
c(t) = \int_{h}^t 2 \frac{a(s)}{b(s)} \, d s + \int_{t-h}^h u(s) \, d s + K
\ee
Evaluating at $t=h$ we see from (\ref{init-b0}) that $K=0$, and hence we have at $t=2h$ the relation
\be
c(2h) = \int_{h}^{2h} 2 \frac{a(s)}{b(s)} \, d s
\ee
Now for $t \ge 2h$ we have
\be
c'(t) = b'(t) - a'(t) = 2 \frac{a(t)}{b(t)} - 2 \frac{a(t-h)}{b(t-h)}
\ee
and thus for some constant $K'$
\be
c(t) = \int_{t-h}^t 2 \frac{a(s)}{b(s)} \, d s + K'
\ee
Evaluating at $t=2h$ we deduce that $K'=0$, and this establishes Property (4). Property (5) follows immediately.

\subsection{Fluid limit: Initial conditions for the tangle from DDE initial condition}\label{sec:init-DDE->stoch}
Let $\alpha$ be a DDE initial condition. As Lemma \ref{lem1} shows, $\alpha$ provides the necessary information
to generate a unique solution of the delay equations (\ref{DDE1}).
We will now show that $\alpha$ also generates the initial conditions for a tangle process.
Recall that $L_m,\dots,L_{2m}$ are determined by the initial conditions
$\xi_m,u_1,\dots,u_m$ through the relation (\ref{init-proc3}),
and that the function $\{b(t) \,:\, h \le t \le 2h\}$ is determined by $a(h), \{u(s):0 \le s \le h\}$ through the formula (\ref{init-b1}).
We will choose the initial values $u_1,\dots,u_m$ for the tangle process depending on the function $u(s)$ in such a way that
the difference $B^{(\lambda)}(t) - b(t)$ is small for
all $t \in [h,2h]$, where $B^{(\lambda)}(t)$ is the 
rescaled variable defined in (\ref{rescale1}). Define the set of all initial value sequences:
\be
{\cal S}(m) = \{v = (v_1,\dots,v_m) \,: \, v_i \in \{0,1,2\}, \, i=1,\dots,m\}
\ee

\begin{defn}
Let $\alpha = (a(h), \{u(s) : 0 \le s \le h\})$ be a DDE initial condition, and let 
$b_{\alpha}(t)$ be defined for $h \le t \le 2h$ by the formula (\ref{init-b1}). 
Let $\xi_{\alpha} = \max (\lfloor \lambda a(h) \rfloor, 1 )$.
Given an initial condition $(\xi_{\alpha}, v)$ for the tangle process, where $v \in {\cal S}(m)$,
let $\{ B_v^{(\lambda)}(t) \,:\, h \le t \le 2h \}$ be given by (\ref{rescale1}) where $L_m,\dots,L_{2m}$ are
defined by the formula (\ref{init-proc3}) with $\xi_m = \xi_{\alpha}$ and $u_i = v_i$. Define
\be\label{def:F}
F(\alpha, \lambda) = \{v \in {\cal S}(m) \,:\, \sup_{m \le n \le 2 m} | B_v^{(\lambda)}(t_n) - b_{\alpha}(t_n) | \le 4 h^{1/2} \lambda^{-1/2} + \lambda^{-1} \}
\ee
\end{defn}

\begin{lemma}\label{lem1a}
The set $F(\alpha, \lambda)$ is non-empty.
\end{lemma}

\medskip
\noindent{\em Proof:}
from (\ref{init-proc3}) and (\ref{init-b1}) we derive for $m \le n \le 2 m$
\be\label{pf:lem2-1}
B_v^{(\lambda)}(t_n) - b_{\alpha}(t_n) &=& \lambda^{-1} \sum_{i=n-m+1}^m v_i - \int_{t_{n-m}}^{t_m} u(s) \, d s
+  \lambda^{-1} \xi_{\alpha} - a(h) \nonumber \\
&=& \lambda^{-1} \sum_{i=n-m+1}^m (v_i - x_i)
+  \lambda^{-1} \xi_{\alpha} - a(h)
\ee
where
\be
x_j = \lambda \, \int_{t_{j-1}}^{t_j} u(s) \, d s, \quad j = 1,\dots, m
\ee
We also have from the definition of $\xi_{\alpha}$
\be\label{pf:lem2-2}
| \lambda^{-1} \xi_{\alpha} - a(h) | \le \lambda^{-1}
\ee
We now introduce a product probability measure on ${\cal S}(m)$ so that the coordinates $v_1,\dots,v_m$ are independent
random variables:
for any sequence $(u_1,\dots,u_m)$,
\be
\P(v = (u_1,\dots,u_m)) = \prod_{j=1}^m \P_j(v_j = u_j)
\ee
The distribution $\P_j$ is chosen so that
\be\label{pf:lem2-3}
\E[v_j] = \sum_{k=0,1,2} k \, \P_j(v_j = k) = x_j
\ee
(since $0 \le x_j \le 2$ this is always possible). Define
\be
M_n = \sum_{j=1}^n (v_j - x_j), \quad 1 \le n \le m
\ee
Since the $\{v_j\}$ are independent with finite variances and (\ref{pf:lem2-3}) holds, we can apply Kolmogorov's maximal inequality \cite{Billingsley}
and deduce that for any $\delta > 0$
\be
\P\left( \max_{1 \le n \le m} |M_n| > \delta\right) \le \delta^{-2} \, {\rm VAR}[M_m]
\ee
Since $|v_j| \le 2$ for all $j$, we have ${\rm VAR}[v_j - x_j] \le 4$, and hence by independence
\be
{\rm VAR}[M_m] \le 4 m = 4 \lambda h
\ee
Taking $\delta = 4 h^{1/2} \, \lambda^{1/2}$ we deduce that
\be
\P\left( \max_{m \le n \le 2 m} | \sum_{i=n-m+1}^m (v_i - x_i)) | > 4 h^{1/2} \, \lambda^{1/2} \right) \le 1/4
\ee
Therefore using (\ref{pf:lem2-2}) and the formula (\ref{pf:lem2-1}) we get
\be
\P(F(\alpha, \lambda)) \ge \P\left( \max_{m \le n \le 2 m} |\lambda^{-1} \sum_{i=n-m+1}^m (v_i - x_i)) | \le 4 h^{1/2} \, \lambda^{1/2} \right) \ge 3/4
\ee
and so we deduce that $F(\alpha, \lambda)$ is non-empty.

\section{Statement of results}\label{sec:results}

\begin{thm}\label{thm1}
Let $\alpha$ be a DDE initial condition, and let $(a_{\alpha}(t),b_{\alpha}(t))$ be the associated solutions of the fluid equations (\ref{DDE1})
as described in Lemma \ref{lem1}.  Let $v \in F(\alpha, \lambda)$, and let
$(A_v^{(\lambda)}(t), B_v^{(\lambda)}(t))$ be the rescaled
tangle process with initial conditions $(\xi_{\alpha},v)$ as described in Sections \ref{sec:gen-init-stoch} and \ref{sec:init-DDE->stoch}.
For all $T \ge 2 h$, and for all $\delta > 0$, there is a constant $C < \infty$ (depending on $T,\alpha$)
and $\lambda_0 < \infty$ (depending on $T, \delta,\alpha$) such that for all $\lambda \ge \lambda_0$
\be\label{thm1:eq1}
\P \left( \sup_{2 h \le t \le T} | B_v^{(\lambda)}(t) - b_{\alpha}(t) | > \delta \right) \le
\P \left( \sup_{h \le t \le T} | A_v^{(\lambda)}(t) - a_{\alpha}(t) | > \delta \right) \le C \, \lambda^{-1} \, \delta^{-2}
\ee
\end{thm}

\medskip
\noindent {\em Remark:} Theorem \ref{thm1} confirms that the rescaled 
processes $(A^{(\lambda)}(t), B^{(\lambda)}(t))$ converge in probability to
the deterministic solutions of the delay equations as $\lambda \rightarrow \infty$. This kind of behavior is familiar for
Markov jump processes. One novelty of Theorem \ref{thm1} is that although the processes are not Markov, due to the
delay time $h$, nevertheless the same kind of limiting behavior holds, albeit with the more complicated delay differential equation.

\medskip
The proof of Theorem \ref{thm1} relies on martingale techniques.
The constants $C$ and $\lambda_0$ that appear in the Theorem depend on $\alpha$, the initial conditions 
for the process.
Simulations of the tangle process  \cite{FKS1} have shown that the delay equations (\ref{DDE1}) give an accurate
representation of the tangle even for relatively small values of $\lambda$.

\medskip
The next result shows that the solution of the delay equation (\ref{DDE1}) converges to a constant
as $t \rightarrow \infty$.

\begin{thm}\label{thm2}
Let $\alpha$ be a DDE initial condition, and let $(a_{\alpha}(t),b_{\alpha}(t))$ be the associated solutions of the fluid equations (\ref{DDE1})
as described in Lemma \ref{lem1}. Define
\be\label{thm2:1}
C_1 &=& \sup_{h \le s \le 2 h} |a_{\alpha}(s) - h|, \nonumber \\
\kappa(u) &=& \max \left\{\frac{3}{4}, \, \exp \left( - \frac{h}{3(u + h)}\right) \right\}, \quad u \ge 0 \nonumber \\
 \mu &=& - \frac{1}{2 h} \log (\kappa(C_1/2))
\ee
Then for all $t \ge 4 h$,
\be\label{thm2:2}
| b(t + h) - 2 h | = | a(t) - h | \le C_1 \, \kappa(C_1/2)^{-3/2} \, e^{ - \mu t}
\ee
\end{thm}

Theorem \ref{thm2} shows that the solutions of the delay equation converge exponentially to their stationary values
with rate at least $\mu$.
This limiting behavior shows that the number of tips behaves as $2 \lambda h$ to leading order  for large arrival rates.

\section{Proof of Theorem \ref{thm1}}\label{sec:pf}
Theorem \ref{thm1} will be proved using standard martingale techniques as presented for example in \cite{DarlingNorris}.
For convenience we will drop the subscripts $v,\, \alpha$ on the variables.
We assume that $\lambda$ is sufficiently large
so that $\xi_{\alpha} = \lambda a(h) \ge 1$. Define
\be
l = \min (h, a(h))
\ee
The quantity $l$ will appear in many of the bounds derived later in this proof, and will represent the effect of the
initial conditions on the constants $C$ and $\lambda_0$ appearing in Theorem \ref{thm1}.
It follows from (\ref{init-proc3}) and (\ref{LXm}) that
\be\label{pf:lower-bound-B}
B^{(\lambda)}(t) \ge l \qquad \mbox{for all $t \ge h$}
\ee
and from (\ref{init-b1}), (\ref{init-b2}) that
\be\label{pf:lower-bound-b}
b(t) \ge l \qquad \mbox{for all $t \ge h$}
\ee
Define for $t \in [h,T]$
\be\label{def:g}
g(t) = \sup_{h \le s \le t} |A^{(\lambda)}(s) - a(s)|
\ee
Note that the quantity of interest in Theorem \ref{thm1} is $\P(g(T) > \delta)$.
We next derive the first inequality in (\ref{thm1:eq1}).
Recall that $b(t) = a(t - h) + h$ for all $t \ge 2h$, and that $L_{n} = X_{n-m} + m$
for all $n \ge 2m$. Therefore if $t \ge 2h$ and $t \in [t_n,t_{n+1})$ we have
\be
B^{(\lambda)}(t) - b(t) = B^{(\lambda)}(t_n) - b(t) = A^{(\lambda)}(t_{n-m}) - a(t-h)
= A^{(\lambda)}(t - h) - a(t-h)
\ee
and therefore
\be\label{thm1:pf3aa}
\sup_{2 h \le s \le t} |B^{(\lambda)}(s) - b(s)| = g(t- h) \le g(t)
\ee
This establishes the first inequality in (\ref{thm1:eq1}), and so reduces the result to
deriving a bound for $\P(g(T) > \delta)$.

\medskip
As a first step we will derive a uniform bound for the difference $B^{(\lambda)}(s) - b(s)$
is terms of the function $g$ and an error term coming from the initial conditions.
Recall that by assumption $v \in F(\alpha, \lambda)$, therefore
\be\label{thm1:pf3a}
\sup_{m \le n \le 2 m} | B^{(\lambda)}(t_n) - b(t_n) | \le 4 h^{1/2} \lambda^{-1/2} + \lambda^{-1}
\ee
Furthermore if $t \in [h,2h]$ and $t \in [t_n,t_{n+1})$ then
\be\label{thm1:pf3b}
| b(t) - b(t_n) | = \left| \int_{t_n}^t (1 - u(s-h)) \, d s \right| \le t - t_n \le \lambda^{-1}
\ee
Therefore (\ref{thm1:pf3a}) and (\ref{thm1:pf3b}) together imply that
\be\label{thm1:pf3c}
\sup_{h \le s \le 2h}  |B^{(\lambda)}(s) - b(s)| \le 4 h^{1/2} \lambda^{-1/2} + 2 \, \lambda^{-1}
\le 6 h^{1/2} \lambda^{-1/2}
\ee
where we have used $ h \lambda \ge 1$. Combining (\ref{thm1:pf3aa}) and (\ref{thm1:pf3c}) we get the uniform bound
\be\label{bound-B}
\sup_{h \le s \le t} | B^{(\lambda)}(s) - b(s) |  \le g(t) + 6 h^{1/2} \lambda^{-1/2} \qquad \mbox{for all $t \ge h$}
\ee

\medskip
Next we will derive a bound for the quantity $A^{(\lambda)}(t) - a(t)$.
For all $j \ge m$ we define
\be
G_{j+1} &=& X_{j+1} - X_j - \E[X_{j+1} - X_j \,|\, {\cal F}(\lambda,j)] \label{def:G} \\
H_{j+1} &=& \E[X_{j+1} - X_j \,|\, {\cal F}(\lambda,j)] - \lambda \, (a(t_{j+1}) - a(t_j)) \label{def:H}
\ee
Then we have
\be\label{pf1:eq2}
\lambda^{-1} \, \sum_{j=m}^{n-1} \left(G_{j+1} + H_{j+1} \right) &=&
\lambda^{-1} \, (X_n - X_m) - (a(t_n) - a(t_m)) \nonumber \\
&=& A^{(\lambda)}(t_n) - a(t_n)
\ee
since $X_m = \lambda a(h) = \lambda a(t_m)$.
The sum $\sum_{j=m}^{n-1} G_{j+1}$ is a martingale (as will be explained below) and we will use 
this fact to bound the probability that it grows too large. The other sum
$\sum_{j=m}^{n-1} H_{j+1}$ is treated as an error term, and will be controlled using coarse bounds.
We will first derive the bounds for this sum, then return to the martingale estimates.

\medskip
From (\ref{def:H}), (\ref{init-proc2}) and (\ref{DDE1}) we get
\be\label{pf1:eq4}
H_{j+1} &=& \E[1 - U_{j+1} \,|\, {\cal F}(\lambda,j)]  - \lambda \, \int_{t_{j}}^{t_{j+1}} \left(1 - 2 \, \frac{a(s)}{b(s)} \right) \, d s \nonumber \\
&=& - 2 \, \frac{X_j}{L_j} + \frac{X_j}{L_j^2}  + 2 \, \lambda \, \int_{t_{j}}^{t_{j+1}}   \frac{a(s)}{b(s)}  \, d s \nonumber  \\
&=& 2 \, \lambda \, \int_{t_{j}}^{t_{j+1}}  \left( \frac{a(s)}{b(s)} -  \frac{A^{(\lambda)}(t_j)}{B^{(\lambda)}(t_j)} \right) \, d s 
+ \lambda^{-1} \frac{A^{(\lambda)}(t_j)}{(B^{(\lambda)}(t_j))^2} \nonumber  \\
&=& 2 \, \lambda \, \int_{t_{j}}^{t_{j+1}}  \left( \frac{a(s)}{b(s)} -  \frac{A^{(\lambda)}(s)}{B^{(\lambda)}(s)} \right) \, d s 
+ \lambda^{-1} \frac{A^{(\lambda)}(t_j)}{(B^{(\lambda)}(t_j))^2} \nonumber  \\
\ee
We write
\be\label{pf1:eq5}
 \frac{a(s)}{b(s)} -  \frac{A^{(\lambda)}(s)}{B^{(\lambda)}(s)} &=&
 \frac{a(s) - A^{(\lambda)}(s)}{b(s)} + \frac{A^{(\lambda)}(s)}{B^{(\lambda)}(s)} \,    \frac{B^{(\lambda)}(s) - b(s)}{b(s)}
 \ee
Using the bounds $A^{(\lambda)}(s) \le B^{(\lambda)}(s)$ and (\ref{pf:lower-bound-B}), (\ref{pf:lower-bound-b}), (\ref{def:g}) and (\ref{bound-B})
we have from (\ref{pf1:eq5}) 
\be
\left| \frac{a(s)}{b(s)} -  \frac{A^{(\lambda)}(s)}{B^{(\lambda)}(s)} \right| \le 2 \, l^{-1} \, g(s) + 6 \,  l^{-1} \, h^{1/2} \, \lambda^{-1/2}
\ee
Therefore we deduce from (\ref{pf1:eq4}) that
\be\label{pf1:eq4a}
| H_{j+1} | \le 4 \, l^{-1} \, g(t_{j+1}) +  12 \, l^{-1} \, h^{ 1/2} \, \lambda^{-1/2} + l^{-1} \, \lambda^{-1} \le 4 \, l^{-1} \, g(t_{j+1}) +  13 \, l^{-1} \, h^{ 1/2} \, \lambda^{-1/2}
\ee
which gives the bound
\be\label{bound-H-sum}
\lambda^{-1} \, \left| \sum_{j=m}^{n-1}  H_{j+1} \right| \le 4 \, l^{-1} \, \lambda^{-1} \,  \sum_{j=m}^{n-1} g(t_{j+1}) +  13 \, l^{-1} \, h^{ 1/2} \, \lambda^{-1/2} \, (t_n - h) \nonumber \\
\le 4 \, l^{-1} \, \lambda^{-1} \,  \sum_{j=m}^{n-1} g(t_{j+1}) +  13 \, l^{-1} \, h^{ 1/2} \, \lambda^{-1/2} \, (T - h)
\ee

\medskip
\noindent Next we use the martingale property to bound the first sum on the left side of (\ref{pf1:eq2}). 
Using (\ref{def:G}) we have
\be
G_{j+1} &=& 1 - U_{+1} - \left(1 - \E[U_{j+1} \,|\, {\cal F}(\lambda,j)] \right) \\
&=& 2 \, \frac{X_j}{L_j} - \frac{X_j}{L_j^2} - U_{j+1}
\ee
It follows that $G_{j+1}$ is ${\cal F}(\lambda,j+1)$-measurable, and
\be
\E[G_{j+1} \,|\, {\cal F}(\lambda,j)] = 0
\ee
Furthermore $|G_{j+1}| \le 2$, so $\{G_j\}$ is a bounded martingale difference series relative to the filtration ${\cal F}(\lambda,n)$.
Therefore $\sum_{j=m}^{n-1} G_{j+1}$ is a martingale, and $\left(\sum_{j=m}^{n-1} G_{j+1}\right)^2$ is a 
submartingale, so we can apply Doob's martingale inequality \cite{Billingsley} to deduce that
for any $N > m$ and any $\theta > 0$
\be
\P\left( \sup_{m \le n \le N} \left|\lambda^{-1} \, \sum_{j=m}^{n-1} G_{j+1} \right| \ge \theta \right) & = & 
\P\left( \sup_{m \le n \le N} \left|  \sum_{j=m}^{n-1} G_{j+1} \right|^2 \ge \lambda^2 \, \theta^2  \right) \nonumber \\
&\le& \lambda^{-2} \, \theta^{-2}  \, \E\left[ \left(\sum_{j=m}^{N-1} G_{j+1} \right)^2\right] \nonumber \\
& = & \lambda^{-2} \, \theta^{-2}  \, \sum_{j=m}^{N-1} \E[ G_{j+1}^2 ] \nonumber \\
& \le & 4 \, \theta^{-2}  \, \lambda^{-2} \, (N-m) \nonumber \\
& = & 4 \, \theta^{-2}  \, \lambda^{-1} \, (T - h)
\ee
Define the event
\be\label{def:E}
E = \left\{\sup_{m \le n \le N} \left|\lambda^{-1} \, \sum_{j=m}^{n-1} G_{j+1} \right| < \theta  \right\}
\ee
so we have
\be
\P(E) \ge 1 - 4 \, \theta^{-2}  \, \lambda^{-1} \, (T - h)
\ee
Combining (\ref{pf1:eq2}), (\ref{bound-H-sum}) and (\ref{def:E}), it follows that
the event $E$ implies that for any $m \le n \le N$,
\be\label{pf1:eq6}
| A^{(\lambda)}(t_n) - a(t_n) | 
& \le & \rho + 4 \, l^{-1} \, \lambda^{-1} \,  \sum_{j=m}^{n-1} g(t_{j+1})
\ee
where
\be\label{pf1:eq7}
\rho =  \theta  + 13 \, l^{-1} \, h^{ 1/2} \, \lambda^{-1/2} \, (T - h)
\ee
Since $g(t) \ge 0$, and (\ref{pf1:eq6}) holds for all $m \le n \le N$, this also implies that
\be
\sup_{m \le k \le n} | A^{(\lambda)}(t_k) - a(t_k) | \le \rho + 4 \, l^{-1} \, \lambda^{-1} \,  \sum_{j=m}^{n-1} g(t_{j+1})
\ee
Furthermore if $t \in [2m,T]$ and $t \in [t_k,t_{k+1})$ we have
\bee
| A^{(\lambda)}(t) - a(t) | &=& | A^{(\lambda)}(t_k) - a(t_k) + a(t_k) - a(t)| \\
& \le & | A^{(\lambda)}(t_k) - a(t_k)| + |a(t_k) - a(t)| \\
& \le & | A^{(\lambda)}(t_k) - a(t_k)| + (t - t_k) \\
& \le & | A^{(\lambda)}(t_k) - a(t_k)| + \lambda^{-1}
\eee
where we used the bound $|a'(s)| \le 1$ for all $s > h$ (which follows from (\ref{DDE1})).
Therefore on the event $E$
\be\label{pf1:eq9}
g(t_n)  \le 
\sup_{m \le k \le n} | A^{(\lambda)}(t_k) - a(t_k) | + \lambda^{-1} 
\le  \rho +  \lambda^{-1}  + 4 \, l^{-1} \, \lambda^{-1} \,  \sum_{j=m}^{n-1} g(t_{j+1})
\ee
Now applying the discrete Gronwall inequality \cite{Agarwal} to (\ref{pf1:eq9})  we deduce that for all $m \le n \le N$,
the event $E$ implies that
\be\label{pf1:eq8}
g(t_n) \le (\rho +  \lambda^{-1}) \, e^{4 \, l^{-1} \, \lambda^{-1}  \, (n - m)} \le ( \rho +  \lambda^{-1}) \, e^{4 \, l^{-1} \, (T-h)}
\ee
Given $\delta > 0$ we choose
\be
\theta &=& \frac{\delta}{3} \, e^{ - 4 \, l^{-1} \, (T-h)} \\
\lambda_0 &=& \max \left\{\theta^{-1}, \left(13 \, \theta^{-1} \, l^{-1} \, h^{1/2} \,(T - h))\right)^2 \right\}
\ee
Then for $\lambda \ge \lambda_0$ we have $\rho +  \lambda^{-1} \le 3 \, \theta$ and
\be
( \rho +  \lambda^{-1}) \, e^{4 \, l^{-1} \, (T-h)} \le \delta
\ee
and hence (\ref{pf1:eq8}) implies that $g(T) \le \delta$ on the event $E$. Therefore
\be
\P\left(g(T) > \delta \right) 
\le 1 - \P(E) \le 4 \, \theta^{-2}  \, \lambda^{-1} \, (T - h)
\ee
and this completes the proof with
\be
C = 36 \, e^{ 8 \, l^{-1} \, (T-h)}  \, (T - h)
\ee

\section{Proof of Theorem \ref{thm2}}\label{sec:pf2}
Recall the delay equation (\ref{DDE1}) for $a(t)$. Applying Lemma \ref{lem1} we get
\be\label{pf2:1}
\frac{d a}{d t} = 1 - 2\, \frac{a(t)}{a(t - h) + h}
\ee
Given the solution $a(t)$ for $t \le T$,  (\ref{pf2:1}) is a linear equation for $a(t)$ in the interval $[T,T+h]$,
and we can write down an explicit solution in terms of the solution in the interval $[T-h,T]$.
Then by translating coordinates the equation (\ref{pf2:1}) can be viewed as providing a map from the space of functions on $[0,h]$
into itself. In order to  prove (\ref{thm2:2}) we will consider instead $a(t) - h$, so define for $t \in [0,h]$
\be\label{pf2:2}
x(t) = \frac{a(T-h+t) - h}{2}, \quad y(t) = \frac{a(T+t) - h}{2}
\ee
then from (\ref{pf2:1}) we derive
\be\label{pf2:3}
\frac{d y}{d t} = \frac{x(t) - 2 y(t)}{2 (x(t) + h)}, \quad y(0) = x(h)
\ee
As explained above, we will view (\ref{pf2:3}) as a map from $x$ to $y$.
Define the functional ${\cal F}$ as the map which takes $x$ to the solution $y$ of the equation (\ref{pf2:3}):
\be\label{pf2:4}
{\cal F}(x)(t) = y(t), \quad 0 \le t \le h
\ee
with the norms
\be\label{pf2:5}
\| x \| = \sup_{0 \le t \le h} |x(t)|, \quad \| y \| = \sup_{0 \le t \le h} | {\cal F}(x)(t) |
\ee
We will prove the following bounds: for all differentiable $x$,
\be\label{pf2:6}
\| {\cal F}(x) \| & \le &  \| x \|  \nonumber \\
\| {\cal F}({\cal F}(x)) \| & \le & \kappa(\| x \|) \, \| x \| 
\ee
where $\kappa$ was defined in (\ref{thm2:1}).
Before proving (\ref{pf2:6}) we note that it implies the bound (\ref{thm2:2}):
indeed for $t \ge 4 h$, there is integer $n \ge 2$ such that $(2 n - 1) h \le t < (2 n+1) h$.
The first inequality in (\ref{pf2:6}) implies that
\be\label{pf2:6a}
\sup_{2 n h \le s \le (2 n + 1) h} \, | a(s) - h | \le \sup_{(2 n - 1) h \le s \le 2 n h} \, | a(s) - h |
\ee
Define $x_0(s) = (a(s) - h)/2$ for $s \in [h,2h]$.
Then  for any $t \in [(2 n - 1) h, (2 n+1) h]$ the inequalities  (\ref{pf2:6}) and (\ref{pf2:6a}) imply
\be
| a(t) - h |  &\le& \sup_{(2 n - 1) h \le s \le 2 n h} \, | a(s) - h | \nonumber \\
&=& 2 \, \| {\cal F}^{\, \circ  (2 n-2)} (x_0) \| \nonumber \\
& \le & 2 \, (\kappa(\| x_0 \|))^{n - 1} \, \| x_0 \| \nonumber \\
& \le &  (\kappa(\| x_0 \|))^{(t-h)/2h - 1} \, C_1  \nonumber \\
&=& e^{ - \mu t} \, C_1 \, \kappa(C_1/2)^{-3/2}
\ee
where we used $C_1 = 2 \| x_0 \|$, and also that $\kappa$ is an increasing function.

\medskip
So we have reduced the proof to (\ref{pf2:6}). Given $x$, let $y$ be the solution of (\ref{pf2:3}),
and let $t \in [0,h]$.
There are three cases:

\medskip
\par\noindent \fbox{Case 1: $y'(t) = 0$}
it follows from (\ref{pf2:3}) that $y(t) = x(t)/2$ and hence
\be
| y(t) | \le \frac{1}{2} \, \| x \|
\ee

\medskip
\par\noindent \fbox{Case 2: $y'(t) > 0$}
Define
\be
S_1 &=& \{ s  \in [0,t) \,:\, y'(s) \le 0 \} \nonumber \\
S_2 &=& \{ s \in (t,h] \, : \, y'(s) \le 0 \} \nonumber \\
t_1 &=& \begin{cases} \sup S_1 & \mbox{if $S_1 \neq \emptyset$} \cr 0 & \mbox{if $S_1 = \emptyset$} \end{cases} \nonumber \\
t_2 &=&  \begin{cases} \inf S_2 & \mbox{if $S_2 \neq \emptyset$} \cr h & \mbox{if $S_2 = \emptyset$} \end{cases}
\ee
Then $y(t_1) < y(t) < y(t_2)$.
By assumption $y'$ is continuous, so if $t_1 > 0$ then $y'(t_1) = 0$, and so $y(t_1) = x(t_1)/2$.
If $t_1=0$ then $y(t_1) = y(0) = x(h)$. Thus in either case
\be
y(t) > y(t_1) \ge \min \{x(t_1)/2, x(h) \}
\ee
Similarly if $t_2 < h$ then $y'(t_2) = 0$, and so $y(t_2) = x(t_2)/2$. If $t_2 = h$ then $y(t_2) = y(h)$
and $y'(h) > 0$, so $y(h) < x(h)/2$. Thus
in either case
\be
y(t) < y(t_2) \le \max \{x(t_2)/2, x(h)/2 \}
\ee
Therefore
\be
| y(t) | \le \max \{x(t_2)/2, x(h)/2, - x(t_1)/2, - x(h) \}
\ee
and so we deduce that
\be
| y(t) | \le \max \left\{ | x(h) |, \, \frac{1}{2} \, \| x \| \right\}
\ee

\medskip
\par\noindent \fbox{Case 3: $y'(t) < 0$}
Define
\be
S_3 &=& \{ s  \in [0,t) \,:\, y'(s) \ge 0 \} \nonumber \\
S_4 &=& \{ s \in (t,h] \, : \, y'(s) \ge 0 \} \nonumber \\
t_3 &=& \begin{cases} \sup S_3 & \mbox{if $S_3 \neq \emptyset$} \cr 0 & \mbox{if $S_3 = \emptyset$} \end{cases} \nonumber \\
t_4 &=&  \begin{cases} \inf S_4 & \mbox{if $S_4 \neq \emptyset$} \cr h & \mbox{if $S_4 = \emptyset$} \end{cases}
\ee
Then $y(t_3) > y(t) > y(t_4)$.
By assumption $y'$ is continuous, so if $t_3 > 0$ then $y'(t_3) = 0$, and so $y(t_3) = x(t_3)/2$.
If $t_3=0$ then $y(t_3) = y(0) = x(h)$. Thus in either case
\be
y(t) < y(t_3) \le \max \{x(t_3)/2, x(h) \}
\ee
Similarly if $t_4 < h$ then $y'(t_4) = 0$, and so $y(t_4) = x(t_4)/2$. If $t_4 = h$ then $y(t_4) = y(h) > x(h)/2$, and thus
in either case
\be
y(t) > y(t_4) \ge \min \{x(t_4)/2, x(h)/2 \}
\ee
Therefore
\be
| y(t) | \le \max \{x(t_3)/2, x(h), - x(t_4)/2, - x(h)/2 \}
\ee
and so we deduce again that for this case
\be
| y(t) | \le \max \left\{ | x(h) |, \, \frac{1}{2} \, \| x \| \right\}
\ee

\medskip
Putting together these three cases we have the bound
\be\label{pf2:8}
| y(t) | \le \max \left\{ | x(h) |, \, \frac{1}{2} \, \| x \| \right\}
\ee
This immediately implies that $\| y \| \le \| x \|$ which is the first inequality in (\ref{pf2:6}).
For the second inequality, we will provide a bound for $|y(h)|$ in terms of $\| x \|$, which will be combined with
(\ref{pf2:8}) to derive (\ref{pf2:6}).
Again we examine several cases.

\medskip
\par\noindent \fbox{Case 4: $y'(h) = 0$}
In this case $y(h) = x(h)/2$ and so $|y(h)| \le \|x\|/2$.

\medskip
\par\noindent \fbox{Case 5: $y'(h) < 0$}
In this case $y(h) > x(h)/2$. We assume that 
$y'(t) < 0$ for all $t \in [0,h]$: if this is not true then as with Case 3 we deduce the existence of $t_3$ such that
$y(h) < y(t_3) = x(t_3)/2$, and  then we have $x(h)/2 < y(h) < x(t_3)/2$, which implies $|y(h)| \le \|x\|/2$.
We also assume that $y(h) > 0$: if $y(h) \le 0$ then the inequality $y(h) > x(h)/2$ implies $|y(h)| \le \| x \|/2$.
Since $y(t)$ is monotone decreasing and $y(h) > 0$ this implies that
$y(t) > 0$ for all $t \in [0,h]$.
Also $y'(t) < 0$ implies
\be
y(t) > \frac{x(t)}{2}, \quad \mbox{and} \quad
y(t) > y(h) > \frac{x(h)}{2} \quad \mbox{for all $t \in [0,h)$}
\ee
Suppose first that there is some  $t \in [0,h)$ such that 
\be
y(t) \le \frac{ 3 x(t)}{4}
\ee
Then
\be
\frac{ 3 x(t)}{4} \ge y(t) > y(h) > \frac{x(h)}{2}
\ee
and therefore
\be
|y(h)| \le \frac{3}{4} \, \| x \|
\ee
If no such $t$ exists then we have
\be
y(t) > \frac{ 3 x(t)}{4} \quad \mbox{for all $t \in [0,h)$}
\ee
and hence (since by assumption $y(t) > 0$)
\be
\frac{d y}{ d t} = - \frac{y - x/2}{x + h} \le - \frac{1}{3} \, \frac{y}{x + h}
\le - \frac{1}{3} \, \frac{y}{\| x \| + h}
\ee
We immediately deduce that
\be
y(h) \le y(0) \, \exp \left( - \frac{h}{3(\| x \| + h)} \right)
\ee
Putting together these two possibilities we get
\be\label{pf2:9}
|y(h)| \le \kappa(\| x \|) \, \| x \|  \quad \mbox{where $\kappa(u) = \max \{\frac{3}{4}, \, \exp ( - \frac{h}{3(u + h)}) \}$}
\ee

\par\noindent \fbox{Case 6: $y'(h) > 0$}
The analysis of this case is identical to Case 5 with some signs reversed, and the same conclusion holds.

\medskip
Combining Cases 4,5,6 we conclude that the bound (\ref{pf2:9}) holds in all cases.
Together with (\ref{pf2:8}) we conclude that
\be\label{pf2:10}
\| {\cal F}(x)(h) \| \le \kappa(\| x \|) \, \| x \| 
\ee
Finally we return to the second inequality in (\ref{pf2:6}), and deduce from (\ref{pf2:10}) that
\bee
\| {\cal F}({\cal F}(x)) \| & \le & \max \left\{ \| {\cal F}(x)(h) \|, \, \frac{1}{2} \| {\cal F}(x) \| \right\} \\
& \le &  \max \left\{ \kappa(\| x \|) \, \| x \|, \, \frac{1}{2} \| {\cal F}(x) \| \right\} \\
& \le & \max \left\{ \kappa(\| x \|) \, \| x \|, \, \frac{1}{2} \| x \| \right\} \\
& = & \kappa(\| x \|) \, \| x \|
\eee

\section{Discussion and future directions}\label{sec:discuss}
Theorem \ref{thm1} confirms that the tangle process converges (in probability) to the solution of the delay differential equation
(\ref{DDE1}).
This convergence was explored using numerical simulations in the paper \cite{FKS1}, and was observed to give an accurate representation of the behavior
even for relatively small values of the arrival rate $\lambda$. There are several interesting questions which arise out of this result.
One question is to describe fluctuations of the rescaled process $A^{(\lambda)}(t)$ around the deterministic solution $a(t)$ of the 
delay differential equation. Theorem \ref{thm1}
shows that the scale of fluctuations is not larger than $\lambda^{-1/2}$. This is also the scale of the central limit theorem,
and it would be interesting to determine if the fluctuations are gaussian in leading order. Another interesting question concerns $h$, the duration of the proof
of work. In this paper we assumed throughout that $h$ is constant, however it would be natural to consider $h$ as a random variable.
Finally the convergence of the tangle model to its fluid limit for other tip selection algorithms is also an interesting problem.

\section*{Acknowledgements}
The author thanks Robert Shorten and Pietro Ferraro for helpful discussions and suggestions.

\end{document}